\documentclass{article}
\newcommand{\be}{\mbox{$\begin{equation}}}
\newcommand{\ee}{\mbox{$\end{equation}}}
\usepackage{amssymb}
\usepackage{amsmath}
\usepackage[dvips]{graphicx}
\usepackage[left=2cm,top=3cm,bottom=3cm,right=2cm]{geometry}

\begin{document}

\begin{center}
{\large\textbf{ON SOME NEW INEQUALITIES INVOLVING GENERALIZED ERD\'ELYI-KOBER FRACTIONAL $q$-INTEGRAL OPERATOR}}
\end{center}

\vspace{4mm}

\begin{center}
{\bf{Daniele Ritelli$^{1,*}$ and Praveen Agarwal$^{2}$}}

\vspace{3mm}
$^{1}$ School of Economics, Management and Statistics\\

Department of Statistics\\
University of Bologna,
via Belle Arti 41\\
 40126 Bologna Italy\\
{\bf{ E-mail: daniele.ritelli@unibo.it}}\\
\vspace{3mm}
{\bf{$^{*}$Corresponding Author}}\\
\vspace{3mm}
$^{2}$Department of Mathematics\\
Anand International College of Engineering\\
 Jaipur-303012,
 India.\\
{\bf{ E-mail: goyal.praveen2011@gmail.com}}\\
\end{center}

\vspace{4mm}
\parindent=0mm
{\bf{Abstract.}} In the present investigation, we aim to establish some inequalities involving generalized Erd$\acute{e}$lyi-Kober fractional $q$-integral operator of the two parameters of deformation $q_{1}$ and $q_{2}$ due to Gaulu$\acute{e}$,
 by following the same lines used by Baleanu and Agarwal \cite {BA} in their recent paper. Relevant connections of the results presented here with those earlier ones are also pointed out.
\\



\medskip

\parindent=4mm

\vskip 2mm
\textbf{1. Introduction and Preliminaries}\\
\vskip 2mm

Throughout this paper,  $\mathbb{N}$, $\mathbb{R}$, $\mathbb{C}$, and  ${\mathbb Z}_0^-$
denote the sets of positive integers, real numbers, complex numbers, and nonpositive integers, respectively,
and $\mathbb{N}_0:= \mathbb{N} \cup \left\{0\right\}$.
The enormous success of the theory of integral inequalities involving various fractional integral operators has stimulated the development
 of a corresponding theory in $q$-fractional integral inequalities (see, \emph{e.g.}, \cite{AN, AN-1, BD, Choi,Z, DMB, D, OO, SU}). In this paper we aim to present some inequalities involving generalized Erd$\acute{e}$lyi-Kober fractional $q$-integral operator of the two parameters of deformation $q_{1}$ and $q_{2}$ due to Gaulu$\acute{e}$ \cite{CG}, by following the same lines used by Baleanu and Agarwal \cite {BA} in their recent paper. Relevant connections of the results presented here with those earlier ones are also pointed out.
\vskip 3mm
  Here, we recall the following definitions (see, \emph{e.g.}, \cite[Section 6]{Sr-Ch-12}) and some earlier works.

\vskip 3mm
The $q$-\emph{shifted factorial} $(a;q)_n$ is defined by
\begin{equation}\tag{1.1}\label{q-shifted factorial}
(a;q)_n:=\left\{ \aligned & 1  \hskip 27 mm  \,\,\,\,\, (n=0) \\
                   & \prod_{k=0}^{n-1}\, \left(1-a\,q^k \right)  \quad  \,\,\,
                  \,\, (n \in  \mathbb N),
                   \endaligned   \right.
\end{equation}
where $a,\,q \in \mathbb C$ and it is assumed that $a \ne q^{-m}$ $\left(m \in {\mathbb N}_0\right)$.

\vskip 3mm
The $q$-{\it shifted factorial} for negative subscript is defined by
\begin{equation}\label{negative q-shifted factorial}\tag{1.2}
(a;q)_{-n}:= \frac{1}{\left(1-a\,q^{-1}\right)\,\left(1-a\,q^{-2}\right) \cdots \left(1-a\,q^{-n}\right)\, } \quad \left(n \in {\mathbb N}_0 \right).
\end{equation}
We also write
 \begin{equation}\label{infinite q-shifted factorial}  \tag{1.3}
 (a;q)_\infty := \prod_{k=0}^\infty\, \left(1-a\,q^k\right) \quad (a,\, q \in \mathbb C;\,\, |q|<1).
\end{equation}
It follows from \eqref{q-shifted factorial},  \eqref{negative q-shifted factorial}  and  \eqref{infinite q-shifted factorial} that
 \begin{equation}\label{finite-infinite} \tag{1.4}
(a;q)_n = \frac{(a;q)_\infty}{ \left(a\,q^n;q \right)_\infty} \quad (n \in \mathbb Z),
\end{equation}
which can be extended to $n=\alpha \in \mathbb C$ as follows:
\begin{equation}\label{arbitrary}  \tag{1.5}
(a;q)_\alpha = \frac{(a;q)_\infty}{ \left(a\,q^\alpha;q \right)_\infty} \quad (\alpha \in \mathbb C;\,\,|q|<1),
\end{equation}
where the principal value of $q^\alpha$ is taken.

\vskip 3mm

We begin by noting that Jackson \cite{Jack-10} was the first to develop $q$-calculus in
a systematic way, this topic is also developed in the recent monograph \cite{Vi}.
\vskip 2mm
   The $q$-\emph{derivative} of a function $f(t)$ is defined by
  \begin{equation}\label{q-derivative}  \tag{1.6}
   D_q\{f(t)\}:=\frac{d_q}{d_q\,t}\{f(t)\} = \frac{f(qt)-f(t)}{(q-1)t}.
   \end{equation}
For $q\to1$ we obtain the definition of differential function, being:
$$ \lim_{q \to 1}\, D_q\{f(t)\} = \frac{d}{dt}\{f(t)\},$$
if we assume that $f(t)$ is a differentiable function.

 The function $F(t)$ is a $q$-\emph{antiderivative} of $f(t)$ if $D_q\{F(t)\}=f(t)$. It is denoted by
 \begin{equation}\label{q-antiderivative}  \tag{1.7}
 \int \, f(t)\, d_q \,t.
 \end{equation}
The  \emph{Jackson integral} of $f(t)$ is thus defined,  \emph{formally}, by
  \begin{equation}\label{Jackson integral}  \tag{1.8}
    \int \, f(t)\, d_q \,t:=(1-q)t\, \sum^{\infty}_{j=0}\,q^j\,f \left(q^j\,t\right),
  \end{equation}
  which can be easily generalized in the Stieltjes sense as follows:
    \begin{equation}\label{generalized Jackson integral}  \tag{1.9}
     \int \, f(t)\, d_q \,g(t)= \sum^{\infty}_{j=0}\,f \left(q^j\,t\right)
     \left( g\left(q^j\,t\right)- g\left(q^{j+1}\,t\right) \right).
    \end{equation}

 Suppose that $0<a<b$. The definite $q$-integral is defined as follows:
  \begin{equation}\label{definite q-integral}   \tag{1.10}
  \int_0^b \, f(t)\, d_q \,t:= (1-q)b\, \sum^{\infty}_{j=0}\,q^j\,f \left(q^j\,b\right)
   \end{equation}
and
\begin{equation}\label{another definite q-integral}  \tag{1.11}
  \int_a^b \, f(t)\, d_q \,t= \int_0^b \, f(t)\, d_q \,t- \int_0^a \, f(t)\, d_q \,t.
   \end{equation}
A more general version of \eqref{definite q-integral}  is given by
\begin{equation}\label{general definite q-integral}   \tag{1.12}
     \int_0^b \, f(t)\, d_q \,g(t)= \sum^{\infty}_{j=0}\,f \left(q^j\,b\right)
     \left( g\left(q^j\,b\right)- g\left(q^{j+1}\,b\right) \right).
    \end{equation}

\vskip 3mm
The classical Gamma function $\Gamma (z)$ (see, \emph{e.g.},  \cite[Section 1.1]{Sr-Ch-12}) was introduced by Leonhard Euler in 1729 while he
was trying to extend the factorial $n!=\Gamma (n+1)$ $\left(n \in {\mathbb N}_0\right)$
to real numbers. The $q$-factorial function $(n)_q!$ $\left(n \in {\mathbb N}_0\right)$ of $n!$   defined by
\begin{equation}\label{q-factorial}   \tag{1.13}
\left(n\right)_q!  :=\left\{ \aligned & 1  \hskip 42 mm \,\, \,\,\,\,\, (n=0), \\
                   & \left(n\right)_q\, \left(n-1\right)_q \cdots  \left(2\right)_q\, \left(1\right)_q\,\quad \,\,\,\,
                          \,\,\,  \,\, (n \in  \mathbb N).
                   \endaligned   \right.
\end{equation}
 can be rewritten as follows:
\begin{equation}\label{alternative-q-factorial}  \tag{1.14}
 (1-q)^{-n}\, \prod_{k=0}^{\infty}\, \frac{\left(1-q^{k+1} \right)}{\left( 1-q^{k+1+n} \right)}
  = \frac{(q;q)_\infty}{\left(q^{n+1};q\right)_\infty}\, (1-q)^{-n}
  := \Gamma_q (n+1) \quad (0<q<1 ).
  \end{equation}
Replacing $n$ by $a-1$ in \eqref{alternative-q-factorial},  Jackson \cite{Jack-10} defined the $q$-Gamma function $\Gamma_q (a)$ by
 \begin{equation}\label{q-Gamma-function}  \tag{1.15}
\Gamma_q (a):= \frac{(q;q)_\infty}{\left(q^a;q\right)_\infty}\, (1-q)^{1-a} \quad (0<q<1).
\end{equation}

\vskip 3mm
  The $q$-analogue of $(t-a)^n$ is defined by the polynomial
\begin{equation}\label{q-analogue-power-function} \tag{1.16}
 \aligned  (t-a)_q^n :=&  \left\{\aligned & 1 \hskip 53 mm \quad  (n=0), \\
               & (t-a)\,(t-q\,a) \cdots (t- q^{n-1}\,a) \hskip 5 mm \quad (n \in \mathbb N).
     \endaligned  \right. \\
     =& t^n\, \left( \frac{a}{t}\,;\,q \right)_n \quad \left(n \in \mathbb{N}_0 \right).
     \endaligned
  \end{equation}
\vskip 3mm
\textbf{Definition 1.}  A real-valued function $f(t)$ $(t>0)$ is said to be in the space $C_\lambda$ $(\lambda\in \mathbb{R})$ if there exists a real number $p>\lambda$ such that $f(t)=t^{p}\,\phi(t)$, where $\phi(t) \in C(0,\,\infty)$.

\vskip 3mm
\textbf{Definition 2.}  A function $f(t)$ $(t>0)$ is said to be in the space $C_\lambda^{n}$ $(n\in \mathbb{R})$ if  $f^{(n)}\in C_\lambda$.

\vskip 3mm

\textbf{Definition 3.}  Let $0<q<1, f\in C_\lambda$. Then for $\Re (\beta), \Re(\mu)>0$,  and $\eta \in \mathbb{C}$ we define a generalized Erd\'elyi-Kober fractional integral $I^{\alpha,\beta,\eta}_{q}$ as follows (see, \cite{CG}) :

\begin{equation}\label{Saigo}\tag{1.17}
\aligned &  I^{\eta, \mu, \beta}_{q}\left\{f(t)\right\}=\frac{\beta\,t^{-\beta(\eta+\mu)}}{\Gamma_{q}(\mu)}
\int_0^t(t^{\beta}-\tau^{\beta}\,q)_{\mu-1}\tau^{\beta(\eta+1)-1}f(\tau)d_{q}\tau,\\
& \hskip 20mm
=\beta(1-q^{1/\beta})(1-q)^{\mu-1}\sum_{k=0}^{\infty}\frac{(q^{\mu};q)_{k}}{(q;q)_{k}}q^{k(\eta+1)}
f(tq^{k/\beta}).
\endaligned
\end{equation}
\vskip 3mm

\textbf{Definition 4.}  Let $0<q<1, f\in C_\lambda$. Then for $\Re(\mu)>0$ and $\eta\in\mathbb{C}$ a $q$-analogue of the Kober fractional integral operator is given by (see \cite{AG})
\begin{equation}\label{q-Erd-Kober}\tag{1.18}
\aligned I^{\eta, \mu}_{q}\left\{f(t)\right\}=\frac{t^{-\eta-\mu}}{\Gamma_{q}(\mu)}
\int_0^t(t-\tau\,q)_{\mu-1}\tau^{\eta}f(\tau)d_{q}\tau.
\endaligned
\end{equation}
\vskip 3mm
\textbf{Remark 1.} It is easy to see that
    \begin{equation}\label{Comparion for q-Saigo-Two} \tag{1.19}
  \Gamma_{q}(\mu)>0; \quad \left(q^{\mu};q\right)_{k} >0,
    \end{equation}
     for all $\mu>0$ and $k \in \mathbb{N}_0$.
    If $f:[0,\infty)\to[0,\infty) $ is a continuous function. Then we conclude that, under the given conditions in \eqref{Saigo}, each term in the series of generalized Erd\'elyi-Kober $q$-integral operator is nonnegative defined by

      \begin{equation}\label{Comparion for q-Saigo} \tag{1.20}
 I^{\eta, \mu, \beta}_{q} \left\{f(t)\right\} \geq 0,
 \end{equation}
 \emph{for all} $\mu>0$ \emph{and} $\eta \in \mathbb{C}$.
 \vskip 3mm
      On the same way each term in the series of Kober $q$-integral operator \eqref{q-Erd-Kober} is also nonnegative
     defined by
     \begin{equation}\label{Comparison-q-Erd-Kober}\tag{1.21}
 I^{\eta, \mu}_{q}\left\{f(t)\right\}\geq 0.
\end{equation}
 \emph{for all} $\mu>0$ \emph{and} $\eta \in \mathbb{C}$.
 \vskip 3mm

 \textbf{2. Generalized Erd\'elyi-Kober $q$-integral Inequalities}\\

\noindent

 In this section, we present six $q$-integral inequalities, which are the core of our research, involving the generalized Erd\'elyi-Kober $q$-integral \eqref{Saigo}
     stated in Theorem 1 to 6 below. Before we recall, as stated in \cite[p. 1, Eq. 2]{Choi} therein, that if $f,\,g$ are two real function defined and integrable on a real interval $\left(I\equiv I\in [a,b]\right),$ we say that $f$ and $g$ are synchronous on $I$ if for each $x,\,y\in I$ the following inequality holds:
\begin{equation}\tag{2.1}\label{I-1}
\left(f(x)-f(y)\right)\left(g(x)-g(y)\right)\geq0 .
\end{equation}

Similarly $f$ and $g$ are asynchronous on $I$ if for any $x,\,y\in I$ the inequality is reversed, that is:
\begin{equation}\tag{2.2}\label{I-2}
\left(f(x)-f(y)\right)\left(g(x)-g(y)\right)\leq0 .
\end{equation}

 {\bf{Theorem 1.}} {\it{Let $0<q_{1}<1$ and $0<q_{2}<1$ with $f,\, g,\, h \in C_{\lambda}$ are three synchronous functions on $[0, \infty)$. If $u : [0,
 \infty)\rightarrow[0, \infty)$ is a continuous function, then the following inequality holds true:

\begin{equation}\tag{2.3}\label{Th-1}
\aligned &
I^{\zeta,\nu, \delta}_{q_{2}} \left\{u(t)\,f(t)\,g(t)\,h(t)\right\}\,I^{\eta, \mu, \beta}_{q_{1}} \left\{u(t)\right\}+
I^{\zeta,\nu, \delta}_{q_{2}} \left\{u(t)\,f(t)\,g(t)\right\}\,I^{\eta, \mu, \beta}_{q_{1}} \left\{u(t)\,h(t)\right\}\\
& \hskip 10mm + I^{\zeta,\nu, \delta}_{q_{2}} \left\{u(t)\,h(t)\right\}\,I^{\eta, \mu, \beta}_{q_{1}} \left\{u(t)\,f(t)\,g(t)\right\}+
I^{\zeta,\nu, \delta}_{q_{2}} \left\{u(t)\right\}\,I^{\eta, \mu, \beta}_{q_{1}} \left\{u(t)\,f(t)\,g(t)\,h(t)\right\}\\
& \hskip 10mm \geq I^{\zeta,\nu, \delta}_{q_{2}} \left\{u(t)\,g(t)\,h(t)\right\}\,I^{\eta, \mu, \beta}_{q_{1}} \left\{u(t)\,f(t)\right\}+
I^{\zeta,\nu, \delta}_{q_{2}} \left\{u(t)\,f(t)\,h(t)\right\}\,I^{\eta, \mu, \beta}_{q_{1}} \left\{u(t)\,g(t)\right\}\\
& \hskip 10mm +
I^{\zeta,\nu, \delta}_{q_{2}} \left\{u(t)\,f(t)\right\}\,I^{\eta, \mu, \beta}_{q_{1}} \left\{u(t)\,g(t)\,h(t)\right\}+
I^{\zeta,\nu, \delta}_{q_{2}} \left\{u(t)\,g(t)\right\}\,I^{\eta, \mu, \beta}_{q_{1}} \left\{u(t)\,f(t)\,h(t)\right\}.
 \endaligned
\end{equation}\\
\emph{for all}  $t>0$, $\mu, \nu, \beta, \delta>0$ \emph{and} $\eta, \zeta \in \mathbb{C}$.}}\\

{\it{Proof}}: Let $f, g$ and $h$ are three continuous and synchronous functions on $[0, \infty)$ and \eqref{I-1} is satisfied. Then for all $\tau, \rho \geq 0$, we have

\begin{equation}\tag{2.4}\label{Th-1.1}
\left(\left(f(\tau)-f(\rho)\right)\left(g(\tau)-g(\rho)\right)\left(h(\tau)+h(\rho)\right)\right)\geq 0,
\end{equation}
which implies that
\begin{equation}\tag{2.5}\label{Th-1.2}
\aligned & f(\tau)g(\tau)h(\tau)+f(\rho)g(\rho)h(\rho)+f(\tau)g(\tau)h(\rho)+f(\rho)g(\rho)h(\tau)\\
 & \geq f(\tau)g(\rho)h(\tau)+f(\tau)g(\rho)h(\rho)+f(\rho)g(\tau)h(\tau)+f(\rho)g(\tau)h(\rho).
 \endaligned
 \end{equation}
Now, multiplying both sides of \eqref{Th-1.2} by
\[
\frac{\beta\,t^{-\beta(\eta+\mu)}}{\Gamma_{q_{1}}(\mu)}
(\tau^{\beta}-\tau^{\beta}\,q_{1})_{\mu-1}\tau^{\beta(\eta+1)-1}u(\tau)
\]
and integrating the resulting inequality with respect to $\tau$ from $0$ to $t$, and using \eqref{Saigo}, we get
\begin{equation}\tag{2.6}\label{Th-1.4}
\aligned &
I^{\eta, \mu, \beta}_{q_{1}} \left\{u(t)\,f(t)\,g(t)\,h(t)\right\}+h(\rho)I^{\eta, \mu, \beta}_{q_{1}}
\left\{u(t)\,f(t)\,g(t)\right\}+f(\rho)g(\rho)I^{\eta, \mu, \beta}_{q_{1}} \left\{u(t)\,h(t)\right\}\\
& \hskip 10mm +f(\rho)g(\rho)h(\rho)I^{\eta, \mu, \beta}_{q_{1}} \left\{u(t)\right\}
 \geq g(\rho)I^{\eta, \mu, \beta}_{q_{1}} \left\{u(t)\,f(t)\,h(t)\right\}+g(\rho)h(\rho)I^{v}_{q_{1}} \left\{u(t)\,f(t)\right\}\\
 & \hskip 60mm +f(\rho)I^{\eta, \mu, \beta}_{q_{1}} \left\{u(t)\,g(t)\,h(t)\right\}+f(\rho)h(\rho)I^{\eta, \mu, \beta}_{q_{1}}
 \left\{u(t)\,g(t)\right\}.
\endaligned
\end{equation}
Next, multiply both sides of \eqref{Th-1.4} by $$\frac{\delta\,t^{-\delta(\zeta+\nu)}}{\Gamma_{q_{2}}(\nu)}
(t^{\delta}-\rho^{\delta}\,q_{2})_{\nu-1}\rho^{\delta(\zeta +1)-1}u(\rho),$$
which remains nonnegative under the conditions in (1.20) and integrating the resulting inequality with respect to $\rho$ from $0$ to $t$, and using \eqref{Saigo}, we
are led to the desired result \eqref{Th-1}. This complete the proof of Theorem 1.\\

{\bf{Theorem 2.}} {\it{Let $0<q_{1}<1,$ $0<q_{2}<1$ and $f, g, h \in C_{\lambda}$ satisfying the condition \eqref{I-1} on $[0, \infty)$. If, we assume $u,\, v : [0,\infty)\rightarrow[0, \infty)$ are continuous functions. Then the following inequality holds true:

\begin{equation}\tag{2.7}\label{Th-2-1}
\aligned &
I^{\zeta,\nu, \delta}_{q_{2}} \left\{v(t)\,f(t)\,g(t)\,h(t)\right\}\,I^{\eta, \mu, \beta}_{q_{1}} \left\{u(t)\right\}+
I^{\zeta,\nu, \delta}_{q_{2}} \left\{v(t)\,f(t)\,g(t)\right\}\,I^{\eta, \mu, \beta}_{q_{1}} \left\{u(t)\,h(t)\right\}\\
& \hskip 10mm + I^{\zeta,\nu, \delta}_{q_{2}} \left\{v(t)\,h(t)\right\}\,I^{\eta, \mu, \beta}_{q_{1}} \left\{u(t)\,f(t)\,g(t)\right\}+
I^{\zeta,\nu, \delta}_{q_{2}} \left\{v(t)\right\}\,I^{\eta, \mu, \beta}_{q_{1}} \left\{u(t)\,f(t)\,g(t)\,h(t)\right\}\\
& \hskip 10mm \geq I^{\zeta,\nu, \delta}_{q_{2}} \left\{v(t)\,g(t)\,h(t)\right\}\,I^{\eta, \mu, \beta}_{q_{1}} \left\{u(t)\,f(t)\right\}+
I^{\zeta,\nu, \delta}_{q_{2}} \left\{v(t)\,f(t)\,h(t)\right\}\,I^{\eta, \mu, \beta}_{q_{1}} \left\{u(t)\,g(t)\right\}\\
& \hskip 10mm +
I^{\zeta,\nu, \delta}_{q_{2}} \left\{v(t)\,f(t)\right\}\,I^{\eta, \mu, \beta}_{q_{1}} \left\{u(t)\,g(t)\,h(t)\right\}+
I^{\zeta,\nu, \delta}_{q_{2}} \left\{v(t)\,g(t)\right\}\,I^{\eta, \mu, \beta}_{q_{1}} \left\{u(t)\,f(t)\,h(t)\right\}.
 \endaligned
\end{equation}\\
\emph{for all}  $t>0$, $\mu, \nu, \beta, \delta>0$ \emph{and} $\eta, \zeta \in \mathbb{C}$.}}\\

{\it{Proof}}: To prove the above result, multiplying both sides of \eqref{Th-1.4} by $$\frac{\delta\,t^{-\delta(\zeta+\nu)}}{\Gamma_{q_{2}}(\nu)}
(t^{\delta}-\rho^{\delta}\,q_{2})_{\nu-1}\rho^{\delta(\zeta +1)-1}v(\rho),$$
which remains nonnegative under the conditions in (1.20) and integrating the resulting inequality with respect to $\rho$ from $0$ to $t$, and using \eqref{Saigo}, we are led to the desired result \eqref{Th-2-1}. This complete the proof of Theorem 2.\\

{\bf{Remark 1.}} It may be noted that the inequalities in \eqref{Th-1} and \eqref{Th-2-1} are reversed if functions $f, g$ and $h$ are asynchronous.
It is also easily seen that the special case $u=v$ of \eqref{Th-2-1} in Theorem 2 reduces to that in Theorem 1.\\
\vskip 3mm
{\bf{Theorem 3.}} {\it{Let $0<q_{1}<1,$ $0<q_{2}<1$ and $u : [0, \infty)\rightarrow[0, \infty)$ be a continuous function with we assume $f, g, h \in C_{\lambda},$ are three synchronous functions on $[0, \infty),$ for which the following condition are satisfied:

\begin{equation}\tag{2.8}\label{Co}
\psi \leq f(x) \leq\Psi,\,\,\,\,\,\phi \leq g(x) \leq \Phi\,\,\, and \,\,\,  \omega \leq h(x) \leq \Omega, \,\,\,\,(\phi, \psi, \omega, \Phi, \Psi,
\Omega \in \mathbb{R}; x \in [0, \infty)).
\end{equation}

Then the following inequality holds true:
\begin{equation}\tag{2.9}\label{Th-3-1}
\aligned &
  \left| ^{\eta, \mu, \beta}_{q_{1}}\left\{u(t)\,f(t)\,g(t)\,h(t)\right\}\,I^{\zeta,\nu, \delta}_{q_{2}}
  \left\{u(t)\right\}+I^{\eta, \mu, \beta}_{q_{1}} \left\{u(t)\,h(t)\right\}\,I^{\zeta,\nu, \delta}_{q_{2}}
  \left\{u(t)\,f(t)\,g(t)\right\}\right.\\
& \left.+I^{\eta, \mu, \beta}_{q_{1}} \left\{u(t)\,g(t)\right\}\,I^{\zeta,\nu, \delta}_{q_{2}} \left\{u(t)\,f(t)\,h(t)\right\}
+I^{\eta, \mu, \beta}_{q_{1}} \left\{u(t)\,f(t)\,\right\}\,I^{\zeta,\nu, \delta}_{q_{2}} \left\{g(t)\,h(t)\,u(t)\right\}\right.\\
& \left.-I^{\eta, \mu, \beta}_{q_{1}} \left\{u(t)\,g(t)\,h(t)\right\}
 I^{\zeta,\nu, \delta}_{q_{2}} \left\{u(t)\,f(t)\right\}
 -I^{\eta, \mu, \beta}_{q_{1}} \left\{u(t)\,f(t)\,h(t)\right\}\,I^{\zeta,\nu, \delta}_{q_{2}} \left\{u(t)\,g(t)\right\}\right.\\
& \left.- I^{\eta, \mu, \beta}_{q_{1}} \left\{u(t)\,f(t)\,g(t)\right\}\,I^{\zeta,\nu, \delta}_{q_{2}} \left\{u(t)\,h(t)\right\}
  -I^{\eta, \mu, \beta}_{q_{1}} \left\{u(t)\right\}\,I^{\zeta,\nu, \delta}_{q_{2}} \left\{u(t)\,f(t)\,g(t)\,h(t)\right\}\right|\\
& \hskip 60mm \leq\,I^{\eta, \mu, \beta}_{q_{1}} \left\{u(t)\right\}\, I^{\zeta,\nu, \delta}_{q_{2}}
\left\{u(t)\right\}\,(\Psi-\psi)(\Phi-\phi)(\Omega-\omega),
 \endaligned
\end{equation}\\
\emph{for all}  $t>0\,\mu,\, \nu,\, \beta,\, \delta>0$ \emph{and} $\eta,\, \zeta \in \mathbb{C}$.}}\\

{\it{Proof}}: Since $f,\, g$ and $h$ are three continuous and synchronous functions on $[0, \infty)$, for all $\tau,\, \rho \geq 0$, the inequality \eqref{I-1} is satisfied. we have from
\eqref{Co}:
$$\left|f(\tau)-f(\rho)\right| \leq (\Psi-\psi), \,\,\,\,\, \left|g(\tau)-g(\rho)\right| \leq (\Phi-\phi),\,\,\, \left|h(\tau)-h(\rho)\right|\leq
(\Omega-\omega),$$
which implies that
\begin{equation}\tag{2.10}\label{Th-3-2}
\left|\left(f(\tau)-f(\rho)\right)\,\left(g(\tau)-g(\rho)\right)\,\left(h(\tau)-h(\rho)\right)\right|\leq  (\Psi-\psi)\,(\Phi-\phi)\,(\Omega-\omega).
\end{equation}
Let us define the function
 \begin{equation}\tag{2.11}\label{Th-3-3}
  \aligned &
  \mathcal{A}(\tau,\rho)=f(\tau)g(\tau)h(\tau)+f(\rho)g(\rho)h(\tau)+f(\tau)g(\rho)h(\rho)+f(\rho)g(\tau)h(\rho)\\
 &- f(\tau)g(\rho)h(\tau)-f(\rho)g(\rho)h(\rho)-f(\tau)g(\tau)h(\rho)-f(\rho)g(\tau)h(\tau),
 \endaligned
 \end{equation}
Multiplying both sides of \eqref{Th-3-3} by
\[
\frac{\beta\,t^{-\beta(\eta+\mu)}}{\Gamma_{q_{1}}(\mu)}
(\tau^{\beta}-\tau^{\beta}\,q_{1})_{\mu-1}\tau^{\beta(\eta+1)-1}u(\tau)
\]
taking $q$-integration of the resulting inequality with respect to $\tau$ from $0$ to $t$ and using \eqref{Saigo}, we
get:
\begin{equation}\tag{2.12}\label{Th-3-4}
\aligned &
\frac{\beta\,t^{-\beta(\eta+\mu)}}{\Gamma_{q_{1}}(\mu)}\int_0^t\,
(\tau^{\beta}-\tau^{\beta}\,q_{1})_{\mu-1}\tau^{\beta(\eta+1)-1}u(\tau)\,\mathcal{A}(\tau,\rho)\,d_{q_{1}}\tau \\
&  = I^{\eta, \mu, \beta}_{q_{1}} \left\{u(t)\,f(t)\,g(t)\,h(t)\right\}+f(\rho)g(\rho)I^{\eta, \mu, \beta}_{q_{1}} \left\{u(t)\,h(t)\right\}
+g(\rho)h(\rho)I^{\eta, \mu, \beta}_{q_{1}}\left\{u(t)\,f(t)\right\}\\
& +f(\rho)h(\rho)I^{\eta, \mu, \beta}_{q_{1}} \left\{u(t)\,g(t)\,\right\}
  -h(\rho)I^{\eta, \mu, \beta}_{q_{1}} \left\{u(t)\,f(t)\,g(t)\right\}- g(\rho)I^{\eta, \mu, \beta}_{q_{1}} \left\{u(t)\,f(t)\,h(t)\right\}\\
 & \hskip 15 mm -f(\rho)I^{\eta, \mu, \beta}_{q_{1}} \left\{u(t)\,g(t)\,h(t)\right\}-f(\rho)g(\rho)h(\rho)I^{\eta, \mu, \beta}_{q_{1}}
 \left\{u(t)\right\}.
                                                                         \endaligned
\end{equation}
Next, multiply both sides of \eqref{Th-3-4} by
\[
\frac{\delta\,t^{-\delta(\zeta+\nu)}}{\Gamma_{q_{2}}(\nu)}
(t^{\delta}-\rho^{\delta}\,q_{2})_{\nu-1}\rho^{\delta(\zeta +1)-1}u(\rho)
\]
then $q$-integrate of the resulting inequality with respect to $\rho$ from $0$ to $t$, and use
\eqref{Saigo}, so that
\begin{equation}\tag{2.13}\label{Th-3-5}
\aligned &
\frac{\beta\,t^{-\beta(\eta+\mu)}}{\Gamma_{q_{1}}(\mu)}\,\frac{\delta\,t^{-\delta(\zeta+\nu)}}{\Gamma_{q_{2}}(\nu)}\int_0^t\,\int_0^t\,
(\tau^{\beta}-\tau^{\beta}\,q_{1})_{\mu-1}\tau^{\beta(\eta+1)-1}u(\tau)\,
(t^{\delta}-\rho^{\delta}\,q_{2})_{\nu-1}\rho^{\delta(\zeta +1)-1}u(\rho)\,\mathcal{A}(\tau,\rho)\,d_{q_{1}}\tau\,d_{q_{2}}\rho\\
&  =I^{\eta, \mu, \beta}_{q_{1}}\left\{u(t)\,f(t)\,g(t)\,h(t)\right\}\,I^{\gamma,\delta,\zeta}_{q_{2}}
\left\{u(t)\right\}+I^{\eta, \mu, \beta}_{q_{1}} \left\{u(t)\,h(t)\right\}\,I^{\gamma,\delta,\zeta}_{q_{2}} \left\{u(t)\,f(t)\,g(t)\right\}\\
& +I^{\eta, \mu, \beta}_{q_{1}} \left\{u(t)\,g(t)\right\}\,I^{\gamma,\delta,\zeta}_{q_{2}} \left\{u(t)\,f(t)\,h(t)\right\}
+I^{\eta, \mu, \beta}_{q_{1}} \left\{u(t)\,f(t)\,\right\}\,I^{\gamma,\delta,\zeta}_{q_{2}} \left\{g(t)\,h(t)\,u(t)\right\}\\
& -I^{\eta, \mu, \beta}_{q_{1}} \left\{u(t)\,g(t)\,h(t)\right\}
 I^{\gamma,\delta,\zeta}_{q_{2}} \left\{u(t)\,f(t)\right\}
 -I^{\eta, \mu, \beta}_{q_{1}} \left\{u(t)\,f(t)\,h(t)\right\}\,I^{\gamma,\delta,\zeta}_{q_{2}} \left\{u(t)\,g(t)\right\}\\
& - I^{\eta, \mu, \beta}_{q_{1}} \left\{u(t)\,f(t)\,g(t)\right\}\,I^{\gamma,\delta,\zeta}_{q_{2}} \left\{u(t)\,h(t)\right\}
  -I^{\eta, \mu, \beta}_{q_{1}} \left\{u(t)\right\}\,I^{\gamma,\delta,\zeta}_{q_{2}} \left\{u(t)\,f(t)\,g(t)\,h(t)\right\}.
                                                            \endaligned
\end{equation}
Finally, by using \eqref{Th-3-2} on to \eqref{Th-3-5}, we arrive at the desired result \eqref{Th-3-1}, involve in Theorem 3, after a little
simplification.\\

{\bf{Theorem 4.}} {\it{Let $0<q_{1}<1,$  $0<q_{2}<1$ and $u,v : [0, \infty)\rightarrow[0, \infty)$ are continuous synchronous functions with $f,\, g,\, h \in C_{\lambda}$ satisfying the condition \eqref{I-1} on $[0, \infty)$ and the following condition:
\begin{equation}\tag{2.14}\label{Th-4}
\psi \leq f(x) \leq\Psi,\,\,\,\,\,\phi \leq g(x) \leq \Phi \, \,\,\,and \,\,\,\, \omega \leq h(x) \leq \Omega \,\,\,\,(\phi, \psi, \omega, \Phi, \Psi,
\Omega \in \mathbb{R}; x \in [0, \infty))
\end{equation}

Then the following inequality holds true:
\begin{equation}\tag{2.15}\label{Th-4-1}
\aligned &
  \left| ^{\eta, \mu, \beta}_{q_{1}}\left\{u(t)\,f(t)\,g(t)\,h(t)\right\}\,I^{\zeta,\nu, \delta}_{q_{2}}
  \left\{v(t)\right\}+I^{\eta, \mu, \beta}_{q_{1}} \left\{u(t)\,h(t)\right\}\,I^{\zeta,\nu, \delta}_{q_{2}}
  \left\{v(t)\,f(t)\,g(t)\right\}\right.\\
& \left.+I^{\eta, \mu, \beta}_{q_{1}} \left\{u(t)\,g(t)\right\}\,I^{\zeta,\nu, \delta}_{q_{2}} \left\{v(t)\,f(t)\,h(t)\right\}
+I^{\eta, \mu, \beta}_{q_{1}} \left\{u(t)\,f(t)\,\right\}\,I^{\zeta,\nu, \delta}_{q_{2}} \left\{g(t)\,h(t)\,v(t)\right\}\right.\\
& \left.-I^{\eta, \mu, \beta}_{q_{1}} \left\{u(t)\,g(t)\,h(t)\right\}
 I^{\zeta,\nu, \delta}_{q_{2}} \left\{v(t)\,f(t)\right\}
 -I^{\eta, \mu, \beta}_{q_{1}} \left\{u(t)\,f(t)\,h(t)\right\}\,I^{\zeta,\nu, \delta}_{q_{2}} \left\{v(t)\,g(t)\right\}\right.\\
& \left.- I^{\eta, \mu, \beta}_{q_{1}} \left\{u(t)\,f(t)\,g(t)\right\}\,I^{\zeta,\nu, \delta}_{q_{2}} \left\{v(t)\,h(t)\right\}
  -I^{\eta, \mu, \beta}_{q_{1}} \left\{u(t)\right\}\,I^{\zeta,\nu, \delta}_{q_{2}} \left\{v(t)\,f(t)\,g(t)\,h(t)\right\}\right|\\
& \hskip 60mm \leq\,I^{\eta, \mu, \beta}_{q_{1}} \left\{u(t)\right\}\, I^{\zeta,\nu, \delta}_{q_{2}}
\left\{v(t)\right\}\,(\Psi-\psi)(\Phi-\phi)(\Omega-\omega),
 \endaligned
\end{equation}\\
\emph{for all}  $t>0,\,\mu,\, \nu,\, \beta,\, \delta>0$ \emph{and} $\eta,\, \zeta \in \mathbb{C}$.}}\\

{\it{Proof}}: Multiplying both sides of \eqref{Th-3-4} by
\[
\frac{\delta\,t^{-\delta(\zeta+\nu)}}{\Gamma_{q_{2}}(\nu)}
(t^{\delta}-\rho^{\delta}\,q_{2})_{\nu-1}\rho^{\delta(\zeta +1)-1}v(\rho)
\]
and taking the
$q$-integration of the resulting inequality with respect to $\rho$ from $0$ to $t$ with the aid of Definition 1 and then applying \eqref{Th-3-2} on
the resulting inequality, we get the desired result \eqref{Th-4-1}.\\

{\bf{Remark 2.}} It is easily seen that the special case $u=v$ of \eqref{Th-4-1} in Theorem 4 reduces to that in Theorem 3.\\

{\bf{Theorem 5.}} {\it{Let $0<q_{1}<1,$ $0<q_{2}<1$ and $u : [0, \infty)\rightarrow[0, \infty)$ be a continuous function. we assume that $f, g,h \in C_{\lambda}$ satisfying the condition \eqref{I-1} on $[0, \infty)$ with constants $L_{1}, L_{2}$ and $L_{3}$, respectively. Then the following inequality holds true:
\begin{equation}\tag{2.16}\label{Th-5-1}
\aligned &
  | I^{\zeta,\nu, \delta}_{q_{2}} \left\{u(t)\,f(t)\,g(t)\,h(t)\right\}\,I^{\eta, \mu, \beta}_{q_{1}} \left\{u(t)\right\}+
  I^{\zeta,\nu, \delta}_{q_{2}} \left\{u(t)\,f(t)\,g(t)\right\}\,I^{\eta, \mu, \beta}_{q_{1}} \left\{u(t)\,h(t)\right\}\\
& \hskip 10mm + I^{\zeta,\nu, \delta}_{q_{2}} \left\{u(t)\,h(t)\right\}\,I^{\eta, \mu, \beta}_{q_{1}} \left\{u(t)\,f(t)\,g(t)\right\}+
I^{\zeta,\nu, \delta}_{q_{2}} \left\{u(t)\right\}\,I^{\eta, \mu, \beta}_{q_{1}} \left\{u(t)\,f(t)\,g(t)\,h(t)\right\}\\
& \hskip 10mm - I^{\zeta,\nu, \delta}_{q_{2}} \left\{u(t)\,g(t)\,h(t)\right\}\,I^{\eta, \mu, \beta}_{q_{1}} \left\{u(t)\,f(t)\right\} -
I^{\zeta,\nu, \delta}_{q_{2}} \left\{u(t)\,f(t)\,h(t)\right\}\,I^{\eta, \mu, \beta}_{q_{1}} \left\{u(t)\,g(t)\right\}\\
& \hskip 10mm -
I^{\zeta,\nu, \delta}_{q_{2}} \left\{u(t)\,f(t)\right\}\,I^{\eta, \mu, \beta}_{q_{1}} \left\{u(t)\,g(t)\,h(t)\right\} -
I^{\zeta,\nu, \delta}_{q_{2}} \left\{u(t)\,g(t)\right\}\,I^{\eta, \mu, \beta}_{q_{1}} \left\{u(t)\,f(t)\,h(t)\right\}|\\
& \hskip 5mm \leq\,L_{1}\, L_{2}\,L_{3}\,\left[I^{\eta, \mu, \beta}_{q_{1}} \left\{\tau^{3}\,u(t)\right\}\, I^{\zeta,\nu, \delta}_{q_{2}}
\left\{u(t)\right\}+3I^{\eta, \mu, \beta}_{q_{1}} \left\{\tau\,u(t)\right\}\, I^{\zeta,\nu, \delta}_{q_{2}} \left\{\rho^{2}u(t)\right\}\right.\\
& \hskip 30mm \left.-3I^{\eta, \mu, \beta}_{q_{1}} \left\{\tau^{2}\,u(t)\right\}\, I^{\zeta,\nu, \delta}_{q_{2}}
\left\{\rho\,u(t)\right\}-I^{\eta, \mu, \beta}_{q_{1}} \left\{u(t)\right\}\, I^{\zeta,\nu, \delta}_{q_{2}} \left\{\rho^{3}\,u(t)\right\}\right].
 \endaligned
\end{equation}\\
\emph{for all}  $t>0$, $\mu, \nu, \beta, \delta>0$ \emph{and} $\eta, \zeta \in \mathbb{C}$.}}\\

{\it{Proof}}: Let us define the following relations for all $\tau, \rho \in [0, \infty) $ :

\begin{equation}\tag{2.17}\label{Th-5-2}
\left|f(\tau)-f(\rho)\right| \leq L_{1}(\tau-\rho), \,\,\,\,\, \left|g(\tau)-g(\rho)\right| \leq L_{2}(\tau-\rho),\,\,\,
\left|h(\tau)-h(\rho)\right|\leq  L_{3}(\tau-\rho),
\end{equation}
which implies that
\begin{equation}\tag{2.18}\label{Th-5-3}
\left|\mathcal{A}(\tau,\rho)\right|\leq L_{1}\,L_{2}\,L_{3}(\tau-\rho)^{3},
\end{equation}
where, $\mathcal{A}(\tau,\rho)$ is given by \eqref{Th-3-3}.
Then, by setting:
\begin{equation}\tag{2.19}\label{Th-5-4}
 \mathcal{B}(\tau,\rho) := L_{1}\,L_{2}\,L_{3}(\tau-\rho)^{3}.
 \end{equation}

First, we multiplying both sides of \eqref{Th-5-4} by
$$\frac{\beta\,t^{-\beta(\eta+\mu)}}{\Gamma_{q_{1}}(\mu)}
(\tau^{\beta}-\tau^{\beta}\,q_{1})_{\mu-1}\tau^{\beta(\eta+1)-1}u(\tau)$$
 and
$$\frac{\delta\,t^{-\delta(\zeta+\nu)}}{\Gamma_{q_{2}}(\nu)}
(t^{\delta}-\rho^{\delta}\,q_{2})_{\nu-1}\rho^{\delta(\zeta +1)-1}u(\rho),$$
respectively, and
taking the $q$-integration of the resulting inequality with respect to $\tau$ and $\rho$ from $0$ to $t$ with the aid of Definition 1 and then
applying \eqref{Th-3-2} and \eqref{Th-5-3} on the resulting inequality, we get the desired result \eqref{Th-5-1}. This complete the proof of Theorem
5.\\

{\bf{Theorem 6.}} {\it{Let $0<q_{1}<1,$ $0<q_{2}<1$ and assume that $u,v : [0, \infty)\rightarrow[0, \infty)$ are continuous functions and suppose $f,\, g,\, h \in C_{\lambda}$ synchronous on $[0, \infty)$ and satisfying \eqref{I-1} with constants $L_{1}, L_{2}$ and $L_{3}$, respectively. Then the following inequality holds true:
\begin{equation}\tag{2.20}\label{Th-6-1}
\aligned &
  | I^{\zeta,\nu, \delta}_{q_{2}} \left\{v(t)\,f(t)\,g(t)\,h(t)\right\}\,I^{\eta, \mu, \beta}_{q_{1}} \left\{u(t)\right\}+
  I^{\zeta,\nu, \delta}_{q_{2}} \left\{v(t)\,f(t)\,g(t)\right\}\,I^{\eta, \mu, \beta}_{q_{1}} \left\{u(t)\,h(t)\right\}\\
& \hskip 10mm + I^{\zeta,\nu, \delta}_{q_{2}} \left\{v(t)\,h(t)\right\}\,I^{\eta, \mu, \beta}_{q_{1}} \left\{u(t)\,f(t)\,g(t)\right\}+
I^{\zeta,\nu, \delta}_{q_{2}} \left\{v(t)\right\}\,I^{\eta, \mu, \beta}_{q_{1}} \left\{u(t)\,f(t)\,g(t)\,h(t)\right\}\\
& \hskip 10mm - I^{\zeta,\nu, \delta}_{q_{2}} \left\{u(t)\,g(t)\,h(t)\right\}\,I^{\eta, \mu, \beta}_{q_{1}} \left\{u(t)\,f(t)\right\} -
I^{\zeta,\nu, \delta}_{q_{2}} \left\{v(t)\,f(t)\,h(t)\right\}\,I^{\eta, \mu, \beta}_{q_{1}} \left\{u(t)\,g(t)\right\}\\
& \hskip 10mm -
I^{\zeta,\nu, \delta}_{q_{2}} \left\{v(t)\,f(t)\right\}\,I^{\eta, \mu, \beta}_{q_{1}} \left\{u(t)\,g(t)\,h(t)\right\} -
I^{\zeta,\nu, \delta}_{q_{2}} \left\{v(t)\,g(t)\right\}\,I^{\eta, \mu, \beta}_{q_{1}} \left\{u(t)\,f(t)\,h(t)\right\}|\\
& \hskip 5mm \leq\,L_{1}\, L_{2}\,L_{3}\,\left[I^{\eta, \mu, \beta}_{q_{1}} \left\{\tau^{3}\,u(t)\right\}\, I^{\zeta,\nu, \delta}_{q_{2}}
\left\{u(t)\right\}+3I^{\eta, \mu, \beta}_{q_{1}} \left\{\tau\,u(t)\right\}\, I^{\zeta,\nu, \delta}_{q_{2}} \left\{\rho^{2}v(t)\right\}\right.\\
& \hskip 30mm \left.-3I^{\eta, \mu, \beta}_{q_{1}} \left\{\tau^{2}\,u(t)\right\}\, I^{\zeta,\nu, \delta}_{q_{2}}
\left\{\rho\,v(t)\right\}-I^{\eta, \mu, \beta}_{q_{1}} \left\{u(t)\right\}\, I^{\zeta,\nu, \delta}_{q_{2}} \left\{\rho^{3}\,v(t)\right\}\right].
 \endaligned
\end{equation}\\
\emph{for all} $t>0$, $\mu, \nu, \beta, \delta>0$ \emph{and} $\eta, \zeta \in \mathbb{C}$.}}\\

{\it{Proof}}: Multiplying both sides of \eqref{Th-5-4} by
$$\frac{\beta\,t^{-\beta(\eta+\mu)}}{\Gamma_{q_{1}}(\mu)}
(\tau^{\beta}-\tau^{\beta}\,q_{1})_{\mu-1}\tau^{\beta(\eta+1)-1}u(\tau)$$
 and
$$\frac{\delta\,t^{-\delta(\zeta+\nu)}}{\Gamma_{q_{2}}(\nu)}
(t^{\delta}-\rho^{\delta}\,q_{2})_{\nu-1}\rho^{\delta(\zeta +1)-1}v(\rho)$$
respectively, and
taking the $q$-integration of the resulting inequality with respect to $\tau$ and $\rho$ from $0$ to $t$ with the aid of Definition 1 and then
applying \eqref{Th-3-2} and \eqref{Th-5-3} on the resulting inequality, we get the desired result \eqref{Th-6-1}. This complete the proof of Theorem
6.\\

{\bf{Remark 3.}} It is easily seen that the special case $u=v$ of \eqref{Th-6-1} in Theorem 6 reduces to that in Theorem 5.\\

 \textbf{Concluding Remarks }\\

We can present a large number of special cases of our main inequalities in Theorems 1 to 6.
Here we  give only one example: Setting $\beta=1$ in Theorem 1- 6, we get the known results given by Baleanu and Agarwal \cite{BA}.
\vskip 3mm

Therefore, we arrived at the conclusion that our present investigation are general in character and useful in deriving various $q$-inequalities in the theory of fractional $q$-integral operators. \\
%

\end{document}